\documentclass[12pt]{article}
\input{epsf.sty}
\usepackage{mathrsfs}
\usepackage{epsfig}
\usepackage{amsmath}
\usepackage{latexsym}
\usepackage{amssymb}
\usepackage{bm}
\usepackage{amsthm}
\theoremstyle{definition}
\newtheorem{theorem}{Theorem}
\newtheorem{lemma}{Lemma}
\newtheorem{assumption}{Assumption}

\newtheorem{corollary}{Corollary}
\newtheorem{definition}{Definition}

\newtheorem{remark}{Remark}
\newtheorem{example}{Example}

\newcommand{\ds}{\displaystyle}
\DeclareMathOperator{\Wra}{Wra} 
\DeclareMathOperator{\Prb}{\mathbb{P}}
\DeclareMathOperator{\Exp}{\mathbb{E}}
\DeclareMathOperator{\diag}{\mathrm{diag}}

\DeclareMathOperator{\sign}{\mathrm{sign}}
\makeatletter
\let\normalequation=\equation
\def\equation{\@ifnextchar[{\subequation}{\normalequation}}
\def\subequation[#1]#2{\@ifundefined{r@#1}%
  {\def\theequation{\bf ??#2}\@warning
    {Reference `#1' on page \thepage \space
     undefined}}{\edef\@tempa{\@nameuse{r@#1}}%
    \edef\theequation{\expandafter\@car\@tempa \@nil#2}}%
  \let\@currentlabel\theequation $$}
\makeatother

\begin{document}

\begin{titlepage}
\begin{center}
{\large\bf Consensus in continuous-time multi-agent systems under
discontinuous nonlinear protocols}
\\[0.2in]
\begin{center}
Liu Bo, Lu Wenlian and Chen Tianping
\end{center}
\begin{center}
School of Mathematical Sciences, Fudan University, Shanghai, 200433,
P.R.China.
\end{center}
\vspace{0.5in}
\end{center}

\begin{abstract}

In this paper, we provide a theoretical analysis for nonlinear
discontinuous consensus protocols in networks of multiagents over weighted
directed graphs. By integrating the analytic tools from nonsmooth stability
analysis and graph theory, we investigate networks with both fixed topology
and randomly switching topology. For networks with a fixed topology, we
provide a sufficient and necessary condition for asymptotic consensus, and
the consensus value can be explicitly calculated. As to networks with
switching topologies, we provide a sufficient condition for the network to
realize consensus almost surely. Particularly, we consider the case that
the switching sequence is independent and identically distributed. As
applications of the theoretical results, we introduce a generalized
blinking model and show that consensus can be realized almost surely under
the proposed protocols. Numerical simulations are also provided to
illustrate the theoretical results.

\bigskip
{\bf Key words:} multiagent systems, consensus, discontinuous, switching,
almost sure
\end{abstract}

\end{titlepage}

\section{Introduction}\hspace{1.25em}
\label{sec:introduction}In many applications involving multiagent systems, groups of
agents are required to agree upon certain quantities of interest.
This is the so-called ``{\em consensus problem}". Due to the broad applications of
multiagent systems, consensus problem arises in various contexts
such as the swarming of honeybees, flocking of birds (Olfati-Saber, 2006),
formation control of autonomous vehicles (Fax \& Murray, 2004),
distributed sensor networks (Cort\'{e}s \& Bullo, 2005) and so on. In the past decades,
a considerable research effort has been devoted to this problem. Various consensus algorithms
have been proposed and studied. For a review, see the survey Olfati-Saber,
Fax \& Murray (2007), Ren, Beard, \& Atkins (2005) and references therein.

Most existing consensus protocols are continuous protocols, i.e.,
the protocol are continuous functions of time $t$ and the states of
the agents. For example, in Olfati-Saber \& Murray (2004), the
authors studied the following linear consensus protocols:
\begin{eqnarray*}
\dot{x}(t)=\sum_{j\in \mathcal{N}_{i}}a_{ij}[x_{j}(t)-x_{i}(t)],
\end{eqnarray*}
where $x_{i}(t)$ is the state of the $i$-th agent at time $t$, and
$\mathcal{N}_{i}$ is the set of neighbors of agent $i$. In Liu,
Chen, \& Lu (2009), the authors studied two types of nonlinear
protocols over directed graphs. The first one is as follows:
\begin{eqnarray}\label{eqnNonlinearProtocol1}
\dot{x}_{i}(t)=\sum_{j=1}^{n}a_{ij}\phi_{ij}(x_{j},x_{i}),~~i=1,2,\cdots,n,
\end{eqnarray}
where $\phi_{ij}$ are nonlinear functions satisfying the following assumption:
\begin{assumption}
\begin{enumerate}
\item $\phi_{ij}$ are locally Lipschitz continuous;
\item $\phi_{ij}(x,y)=0$ if and only if $x=y$;
\item $(x-y)\phi_{ij}(x,y)<0$, $\forall$ $x\ne y$.
\end{enumerate}
\end{assumption}
They prove that this protocol can realize consensus if and only if
the underlying graph has a spanning tree. The second one is as
follows:
\begin{eqnarray}\label{eqnNonlinearProtocol2}
\dot{x}_{i}(t)=-\sum_{j=1}^{n}l_{ij}[h(x_{j})-h(x_{i})],
\end{eqnarray}
where $h$ is a strictly increasing nonlinear function, and the Laplacian matrix $L=[l_{ij}]$ has the form
\begin{eqnarray*}
\left[
\begin{array}{cc}
L_{11} & 0 \\
L_{21} & L_{22}
\end{array}
\right],
\end{eqnarray*}
where $L_{11}$, $L_{22}$ is irreducible, and $L_{21}\ne 0$. They
prove that this protocol can realize consensus value which is a
convex combination of component states of the initial value.

Previous protocols are for static networks, i.e., networks with
fixed topologies. Yet many real world networks are not static. For
example, in a network of mobile agents, the topology of the network
is dynamical due to limited transmission range and the movement of
the agents. In some cases, the network topology changes gradually.
In other cases, it changes abruptly, which induces discontinuity in
the network topology.
%In Moreau (2004), the author studies the following linear consensus protocols in time-varying systems:
%\begin{eqnarray*}
%\dot{x}(t)=A(t)x(t),
%\end{eqnarray*}
%where $A(t)$ is a Metzler matrix with zero row sums and piecewise continuous in $t$. It was proved
%that if there exist a node $k$, a threshold value $\delta>0$, and an interval length $T>0$ such that for all
%$t\in\mathbb{R}$, the $\delta$-graph associated to $\int_{t}^{t+T}A(s)ds$ has the property that all nodes can be reached from node $k$, then consensus can be reached exponentially.

An important class of discontinuous dynamical network topology is
the so-called {\em switching topology}. Let
$0=t_{0}<t_{1}<\cdots<t_{k}<t_{k+1}<\cdots$ be a partition of
$[0,+\infty)$, on each time interval $[t_{k},t_{k+1})$, the network
has a fixed topology, while at each time point $t_{k}$, the topology
switches to another one randomly or according to some given rule.
Linear consensus protocols over networks with stochastically
switching topologies such as independent and identically distributed
switching (Salehi \& Jadbabaie, 2007), Markovian switching (Matei,
Martins, \& Baras, 2008), and adapted stochastic switching (Liu, Lu,
\& Chen, 2011) have been studied and  conditions for almost sure
consensus have been obtained, which indicates that a directed
spanning tree in the expectation is sufficient for almost sure
consensus.

The above mentioned discontinuous consensus protocols are
discontinuous in time $t$ and continuous in the states of the
agents. Besides, there are another important class of discontinuous
consensus protocols which are discontinuous in the states of the
agents, too. Recently, such protocols have been discussed in several
papers. In Cort\'{e}s (2006), based on normalized and signed
gradient dynamical systems associated with the Laplacian potential,
the author proposed the following two discontinuous consensus
protocols:
\begin{eqnarray}
\dot{p}_{i}(t)=\frac{\sum_{j\in \mathcal{N}_{i}}(p_{j}(t)-p_{i}(t))}{\|LP(t)\|_{2}},\label{protocolCenteral}\\
\dot{p}_{i}(t)=\sign\bigg(\sum_{j\in \mathcal{N}_{i}}(p_{j}(t)-p_{i}(t))\bigg)\label{protocolMaxMin},
\end{eqnarray}
where $L$ is the graph Laplacian of the underlying graph, and
$P(t)=[p_{1}(t),\cdots,p_{n}(t)]^{\top}$. Finite time convergence of
both protocols on connected undirected graphs was proved, where the
centralized protocol \eqref{protocolCenteral} can realize average
consensus, while the distributed algorithm \eqref{protocolMaxMin}
can reach average-max-min consensus. In Cort\'{e}s (2008), the
author further considered the following two discontinuous protocols:
\begin{eqnarray}
\dot{p}_{i}=\sign_{+}\bigg(\sum_{j=1}^{n}a_{ij}(p_{j}-p_{i})\bigg),\label{protocolMax}\\
\dot{p}_{i}=\sign_{-}\bigg(\sum_{j=1}^{n}a_{ij}(p_{j}-p_{i})\bigg),\label{protocolMin}
\end{eqnarray}
where $\sign_{+}(x)=0$ if $x\le 0$ and $\sign_{+}(x)=1$ if $x>0$,
$\sign_{-}(x)=0$ if $x\ge 0$ and $\sign_{-}(x)=-1$ if $x<0$. Both
protocols can realize finite time consensus in a strongly connected
weighted directed graph, where protocol \eqref{protocolMax} can
reach max consensus, while protocol \eqref{protocolMin} can reach
min consensus. In Hui, et al. (2008), the author studied the
stability of consensus under the following discontinuous protocol:
\begin{eqnarray*}
\dot{x}_{i}(t)=\sum_{j=1}^{q}C_{(i,j)}\sign(x_{j}-x_{i}).
\end{eqnarray*}
Under the assumption that $C$ is symmetric and
$\mathrm{rank}(C)=q-1$, they proved finite time convergence for this
protocol.

In this paper, we investigate a new type of nonlinear discontinuous protocols, which can be formulated as follows:
\begin{eqnarray*}
\dot{x}_{i}=-\sum_{j=1}^{n}l_{ij}[g(x_{j})-g(x_{i})],~~i=1,\cdots,n,
\end{eqnarray*}
where $L=[l_{ij}]$ is the underlying graph Laplacian, and $g(\cdot)$
is a discontinuous function that will be specified later. First, we
consider networks with fixed topology. Compared to existing works
which only consider connected undirected graphs or strongly
connected directed graphs, we consider more general directed graphs
that has spanning trees. We show that a directed spanning tree is
sufficient for the network to realize asymptotic consensus. And this
condition is not only sufficient but also necessary. This is an
important improvement since directional communication is important
in practical applications and can be easily incorporated, for
example, via broadcasting. Moreover, a lot of important real world
networks such as the leader-follower networks are not strongly
connected. Then, motivated by the work in synchronization analysis
by Lu and Chen (2004),  we locate the consensus value based on the
left eigenvector corresponding to the zero eigenvalue of the graph
Laplacian. Finally, we show that if the consensus value is a
discontinuous point of $g$, and the underlying graph is strongly
connected, then finite time convergence can be realized.

We also consider the consensus protocol over networks with switching
topologies. The time interval between each successive switching is
assumed to be an independent and identically distributed random
variable. And the network topology is also a random sequence. We
prove a sufficient condition for the network to achieve consensus
almost surely in terms of the scramblingness of the underlying
graph. Based on this result, we study the special case where the
switching sequence is independent and identically distributed. We
show that if the underlying graph has a positive probability to be
scrambling, then the protocol can realize consensus almost surely.
Our results indicate that for a network with stochastically
switching topology to reach consensus almost surely, the network is
unnecessary  connected at each time point. This is more general than
the work in Hui, et. al.(2008) on network with switching topology.

Finally, as applications of the theoretical results. We study
consensus in a general blinking network model under the proposed
consensus protocol. Numerical simulations are also provided to
illustrate the theoretical results.

This paper is organized as follows. In Section
\ref{secPreliminaries}, some preliminary definitions and lemmas
concerning graph theory, matrix theory nonsmooth analysis, and
probability, are provided. Consensus analysis under nonlinear
discontinuous protocols with both fixed topology and switching
topology, are carried out in Section \ref{secMainResults}. An
application of the theoretical results to a general blinking network
model with numerical simulations are given in Section
\ref{secNumericalSimulations}. The paper is concluded in Section
\ref{secConclusions}.

\section{Preliminaries}\hspace{1.25em}
\label{secPreliminaries}In this section, we present some definitions
and basic lemmas that will be used later.

\subsection{Algebraic graph theory and matrix theory}

A {\em weighted directed graph} of order $n$ is denoted by a triple
$\{\mathcal{V},\mathcal{E},W\}$ where
$\mathcal{V}=\{v_{1},\cdots,v_{n}\}$ is the vertex set and
$\mathcal{E}\subseteq \mathcal{V}\times \mathcal{V}$ is the edge
set, i.e., $e_{ij}=(v_{i},v_{j})\in \mathcal{E}$ if there is an edge
from $v_{i}$ to $v_{j}$, and $W=[w_{ij}]$, $i,j=1,\cdots, n$, is the weight
matrix which is a nonnegative matrix such that for
$i,j\in\{1,\cdots,n\}$, $w_{ij}>0$ if and only if $i\ne j$ and
$e_{ji}\in\mathcal{E}$. For a weighted directed graph $\mathcal{G}$ of order $n$, the graph
Laplacian $L(\mathcal{G})=[l_{ij}]_{i,j=1}^{n}$ can be defined from the weight
matrix $W$ in the following way:
\begin{eqnarray*}
l_{ij}=\left\{
\begin{array}{cc}
-w_{ij} & i\ne j\\
\sum\limits_{j=1,\atop{j\ne i}}^{n}w_{ij}& j=i.
\end{array}
\right.
\end{eqnarray*}
And for a given Laplacian matrix $L$, the weighted directed graph corresponding to $L$ is written as $\mathcal{G}(L)$.

In this paper, we only consider simple
graphes, i.e., there are no self links and multiple edges. A
directed path of length $r$ from $v_{i}$ to $v_{j}$ is an ordered
sequence of $r+1$ distinct vertices $v_{k_{1}},\cdots, v_{k_{r+1}}$
with $v_{k_{1}}=v_{i}$ and $v_{k_{r+1}}=v_{j}$ such that
$(v_{i_{s}},v_{i_{s+1}})\in \mathcal{E}$. A ({\em directed}) {\em
spanning tree} is a directed graph such that there exists a vertex
$v_{r}$, called the root vertex, such that for any other vertex
$v_{i}\in \mathcal{V}$, there exists a directed path from $v_{r}$ to
$v_{i}$. We say a graph $\mathcal{G}$ has a spanning tree if a
subgraph of $\mathcal{G}$ that has the same vertex set with
$\mathcal{G}$ is a spanning tree. A graph $\mathcal{G}$ is strongly connected
if for any pair of vertices, say, $v_{i}$, $v_{j}$, there exist directed
paths both from $v_{i}$ to $v_{j}$ and from $v_{j}$ to $v_{i}$.

If a graph has spanning trees, then the vertices of the graph can be divided
into two disjoint sets: $S_{1}$, $S_{2}$, where $S_{1}$ contains the vertices
that can be the root of some spanning tree, $S_{2}$ contains all other vertices.
We have the following lemma.
\begin{lemma}\label{lemSpanningTree}
If a graph $\mathcal{G}$ of $n$ vertices has spanning trees, let $S_{1}$, $S_{2}$ be defined
as above, then
\begin{itemize}
\item[(i)]The subgraph of $\mathcal{G}$ induced by $S_{1}$ is strongly connected.
\item[(ii)] $\mathcal{G}$ is strongly connected if and only if $\#S_{1}=n$.
\end{itemize}
\end{lemma}
{\bf \noindent Proof:~}
(i): First, for any given vertices $v_{1},v_{2}\in S_{1}$,
since $v_{1}$ can be the root of some spanning tree, then from definition,
there is a directed path from $v_{1}$ to $v_{2}$. On the other hand, $v_{2}$
can also be the root of some spanning tree, so there also exists a directed
path from $v_{2}$ to $v_{1}$. Second, we prove that these two paths
contain no vertices outside $S_{1}$. Otherwise, there exists a vertex
$v_{3}\not \in S_{1}$ such that $v_{3}$ is on one of the paths. Suppose
$v_{3}$ is on the path from $v_{1}$ to $v_{2}$, then there
is a directed path from $v_{3}$ to $v_{2}$. Since $v_{2}$ is a root, there
exist directed paths from $v_{2}$ to all other vertices. Thus
there are directed paths from $v_{3}$ to all other vertices, which implies
$v_{3}$ also can be the root of some spanning tree. This contradicts the fact
that $v_{3}\not \in S_{1}$.

(ii): If $\#S_{1}=n$, then the subgraph induced by $S_{1}$ is $\mathcal{G}$
itself. From (i), $\mathcal{G}$ is strongly connected. On the other hand,
if $\mathcal{G}$ is strongly connected, from definition, each vertex can
be the root of some spanning tree. Thus, $\#S_{1}=n$.
$\square$

\begin{remark}
It is known that in a leader-follower system, only the leader can
influence the follower, but the follower can not influence the
leader. So the final state of the system is determined only by the
leader. In a strongly connected system, each agent can be seen as a
leader. So the final state of the system is determined by all
agents. Yet there are also many intermediate cases between these two
extremes. In such cases, there are group of leaders, but the whole
system is not strongly connected. Lemma \ref{lemSpanningTree}
unifies these three cases into a general framework.
\end{remark}

From the proof of Lemma \ref{lemSpanningTree}, we can see that there
exist no edges from vertices of $S_{2}$ to vertices of $S_{1}$. then
after a proper renumbering of its vertices, the graph Laplacian $L$
of $\mathcal{G}$ can be written in the following form:
\begin{eqnarray}
L=\left[
\begin{array}{cc}
L_{1}&0\\
*& L_{2},
\end{array}
\right]\label{spanningtree}
\end{eqnarray}
where the square submatrix $L_{1}$ corresponds to the vertex set
$S_{1}$. Since the subgraph induced by $S_{1}$ is strongly
connected, $L_{1}$ is irreducible. By Perron-Frobenius theory, the
left eigenvector of $L_{1}$ corresponding to the eigenvalue $0$ is
positive. Thus we can define the following
\begin{definition}(Weighted root average)\label{defWeightedRootAverage}
Let $L=[l_{ij}]_{i,j=1}^{n}$ be the graph Laplacian of some weighted
directed graph $\mathcal{G}(L)$. Suppose that $\mathcal{G}$ has spanning trees
and $L$ is of the form (\ref{spanningtree}). Let $\xi=[\xi_{1},\cdots,
\xi_{n_{1}}]^{\top}$ be the positive left eigenvector of
$L_{1}$ corresponding to the eigenvalue $0$ such that
$\sum_{i=1}^{n_{1}}\xi_{i}=1$, where $n_{1}=\#S_{1}$. Given some
$x=[x_{1},\cdots,x_{n}]^{\top}\in \mathbb{R}^{n}$, the {\bf weighted
root average} of $x$ with respect to $L$ is defined as:
\begin{equation*}
\Wra(x,L)=\sum_{i=1}^{n_{1}}\xi_{i}x_{i}.
\end{equation*}
\end{definition}

\begin{remark}
In a leader-follower system, the final state of the system is determined by
the leader only. In the case that there are group of leaders, the final state
of the system is determined by the leader group. The weighted root average is
also a generalization from the case of one leader to the case of leader group.
\end{remark}

\begin{figure}
\centering
\includegraphics[width=0.5\textwidth]{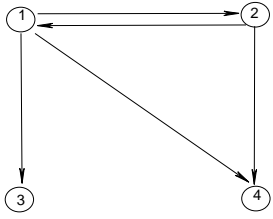}
\caption{Graph example 1}
\label{figGraphExample_1}
\end{figure}

\begin{example}
Consider the graph in Fig. \ref{figGraphExample_1}, it is obvious that
$S_{1}=\{v_{1},v_{2}\}$, $S_{2}=\{v_{3},v_{4}\}$. If we take all the positive
weight of the edges to be $1$, then the graph Laplacian is
\begin{eqnarray*}
L=\left[
\begin{array}{cccc}
1 & -1 & 0 & 0 \\
-1 & 1 & 0 & 0 \\
-1 & 0 & 1 & 0 \\
-1 & -1 & 0 & 2
\end{array}
\right].
\end{eqnarray*}
Here, $L_{1}=\left[
\begin{array}{cc}
1 & -1 \\
-1 & 1
\end{array}\right]$, and $\xi=[1/2,1/2]^{\top}$. Thus, for any
$x=[x_{1},x_{2},x_{3},x_{4}]^{\top}$, $\Wra(x,L)=(x_{1}+x_{2})/2$.
\end{example}
A {\bf Metzler matrix} is a matrix that has nonnegative off-diagonal
entries. It is clear that $-L$ is a Metzler
matrix with zero row sum. Following Liu and Chen (2008), for a Metzler matrix
$M=[m_{ij}]$, we define a
function %$\eta(\cdot)$ as follows:
\begin{eqnarray*}
\eta(M)=\max_{i,j}\{-(m_{ij}+m_{ji})-\sum_{k\ne
i,j}\min\{m_{ik},m_{jk}\}\}.
\end{eqnarray*}
and we say that $M$ is {\bf scrambling} if $\eta(M)<0$. It is
obvious that scramblingness is not influenced by the diagonal
entries of a Metzler matrix, so $L$  is scrambling if and only if
$W$ is scrambling. On the other hand, since there is a one to one correspondence
between each weighted directed graph and its weight matrix $W$ (or Laplacian
matix $L$), we also say a graph is scrambling if $W$ (or $-L$) is scrambling.

\begin{remark}
It can be seen from the definition that if a graph is scrambling,
then for each vertex pair $(v_{i},v_{j})$, either there exists at
least one directed edge between $v_{i}$ and $v_{j}$, or there is
another vertex $v_{k}$ such that there are directed edges from
$v_{k}$ to $v_{i}$ and $v_{k}$ to $v_{j}$. From this, it can be seen
that the graph in Fig. \ref{figGraphExample_1} is scrambling. Since
there exists directed edges between $(v_{1},v_{2})$,
$(v_{1},v_{3})$, $(v_{1},v_{4})$, and $(v_{2},v_{4})$. And there
exist edges from $v_{1}$ to $v_{2}$, $v_{3}$, and edges from $v_{1}$
to $v_{3}$, $v_{4}$.
\end{remark}

If we incorporate a positive threshold $\delta$ on the graph
$\mathcal{G}$, then we get the concept of {\bf $\delta$-graph}
(Moreau, 2004). The $\delta$-graph of $\mathcal{G}$ is a graph that
has the same vertex set and weight matrix with $\mathcal{G}$. Yet
for each $v_{i}$, $v_{j}$, there is a directed edge from $v_{j}$ to
$v_{i}$ if and only if $w_{ij}\ge \delta$. We say a graph
$\mathcal{G}$ is {\bf $\delta$-scrambling} if its $\delta$-graph is
scrambling.
\begin{remark}\label{remScrambling}
It is obvious that if $\mathcal{G}$ is $\delta$-scrambling, then $\eta(-L(\mathcal{G}))\ge\delta$.
\end{remark}
\subsection{Nonsmooth stability analysis}\label{secNonsmooth}
In this subsection, we will provide some concepts and lemmas
concerning nonsmooth stability analysis. First, we present
some basic concepts and theorems from Filippov theory on
differential equations with discontinuous righthand sides.
For more details, the readers are referred to Filippov (1988) directly.

Consider the following differential equations:
\begin{eqnarray}\label{eqnFilippovDefinition}
\dot{x}(t)=f(x(t))
\end{eqnarray}
where $x\in {\mathbb R}^{n}$, and $f$: ${\mathbb R}^{n}\mapsto
{\mathbb R}^{n}$ is a discontinuous map. Then the Filippov solution
of (\ref{eqnFilippovDefinition}) can be defined as:
\begin{definition}\label{defFilippovSolution}
An absolutely continuous function $\varphi$: $[t_{0},t_{0}+a]\mapsto
{\mathbb R}^{n}$ is said to be a Filippov solution to
(\ref{eqnFilippovDefinition}) on $[t_{0}, t_{0}+a]$ if it is a
solution of the differential inclusion:
\begin{eqnarray}\label{eqnDifferentialInclusion}
\dot{x}(t) \in
\bigcap_{\delta>0}\bigcap_{\mu(N)=0}\mathcal{K}[f(B(x,\delta)\backslash
N)],\qquad a.e. t\in [t_{0},t_{0}+a],
\end{eqnarray}
where $\mathcal{K}(E)$ is the closure of the convex hull of $E$,
$B(x,\delta)$ is the open ball centered at $x$ with radius
$\delta>0$, and $\mu(\cdot)$ denote the usual Lebesgue measure in
$\mathbb{R}^{n}$.
\end{definition}
For the simplicity of notation, we denote
$\mathcal{K}[f](x)=\bigcap_{\delta>0}\bigcap_{\mu(N)=0}\mathcal{K}[f(B(x,\delta)\backslash
N)]$, and (\ref{eqnDifferentialInclusion}) can be rewritten as:
\begin{eqnarray}\label{eqnDifferentialInclusionNewVersion}
\dot{x}(t)\in \mathcal{K}[f](x(t)),\qquad a.e. t\in [t_{0},t_{0}+a].
\end{eqnarray}
A Filippov solution of (\ref{eqnDifferentialInclusionNewVersion}) is
a maximum solution if its domain of existence is maximum, i.e., it
can not be extended any further. A set $S\subseteq \mathbb{R}^{n}$
is weakly invariant (resp. strongly invariant) with respect to
(\ref{eqnDifferentialInclusionNewVersion}) if for each $x_{0}\in S$,
$S$ contains a maximum solution (resp. all maximum solutions) from
$x_{0}$ of (\ref{eqnDifferentialInclusionNewVersion}).

Let $f$: $\mathbb{R}^{n}\mapsto \mathbb{R}$, then the usual
one-sided directional derivative of $f$ at $x$ in direction $v$ is
defined as:
\begin{eqnarray}
f'(x,v)=\lim_{t\to 0^{+}}\frac{f(x+tv)-f(x)}{t}.
\end{eqnarray}
The generalized directional derivative of $f$ at $x$ in direction
$v$ is defined as:
\begin{eqnarray}
f^{\circ}(x,v)=\limsup_{y\to x,t\to 0^{+}}\frac{f(y+tv)-f(y)}{t}.
\end{eqnarray}
\begin{definition}(Clarke,1983)\label{defRegularity}
Let $f$: $\mathbb{R}^{n}\mapsto \mathbb{R}$, $f$ is said to be
regular at $x$ if for all $v\in\mathbb{R}^{n}$, the usual one-sided
directional derivative $f'(x,v)$ exists, and $f'(x,v)=f^{
\circ}(x,v)$.
\end{definition}

Following lemma can be used to derive regularity.
\begin{lemma}(Clarke,1983)\label{lemJudgeRegularity}
Let $f$: $\mathbb{R}^{n}\mapsto \mathbb{R}$ be Lipschitz near $x$,
then
\begin{enumerate}
\item If $f$ is convex, then $f$ is regular at $x$;
\item A finite linear combination (by nonnegative scalars) of
functions regular at $x$ is regular at $x$.
\end{enumerate}

\end{lemma}

From Rademacher's Theorem (Clarke,1983), we know that locally
Lipschitz functions are differentiable almost everywhere.

\begin{definition}(Clarke,1983)\label{defGeneralizedGradientLieDerivative}
Let $V$: $\mathbb{R}^{n}\mapsto \mathbb{R}$ be a locally Lipschitz
continuous function. Let $\Omega_{V}$ be the set of points where $V$
fails to be differentiable, then the {\bf Clarke generalized
gradient} of $V(x)$ at $x$ is the set
\begin{eqnarray}
\partial V(x)\triangleq \{\lim_{i\to+\infty}\nabla V(x^{i}): x^{i}\to x, x^{i}\not\in \Omega_{V}\cup
\mathcal{S}\}
\end{eqnarray}
where $\mathcal{S}$ can be any set of zero measure. The {\bf
set-valued Lie derivative } of $V$ with respect to
(\ref{eqnDifferentialInclusionNewVersion}) at $x$ is:
\begin{eqnarray}
\tilde{\mathcal{L}}_{f}V(x)=\{a\in\mathbb{R}: \exists v\in
\mathcal{K}[f](x) \textnormal{ such that } a=\zeta\cdot v,~~ \forall
\zeta\in\partial V(x)\}.
\end{eqnarray}
\end{definition}

The following lemma shows that the evolution of the Filippov
solutions can be measured by the Lie derivative.

\begin{lemma}\label{lemEvolutionAlongFilippovSolution}
Let $x$: $[t_{0},t_{1}]$ be a Filippov solution of
(\ref{eqnFilippovDefinition}). Let $V$: $\mathbb{R}^{n}\mapsto
\mathbb{R}$ be a locally Lipschitz and regular function. Then,
$t\mapsto V(x(t))$ is absolutely continuous,
$\ds\frac{dV(x(t))}{dt}$ exists a.e. and
$\ds\frac{dV(x(t))}{dt}\in\tilde{\mathcal{L}}_{f}V(x(t))$ for a.e.
$t$.
\end{lemma}

\begin{figure}
\centering
\includegraphics[width=0.5\textwidth]{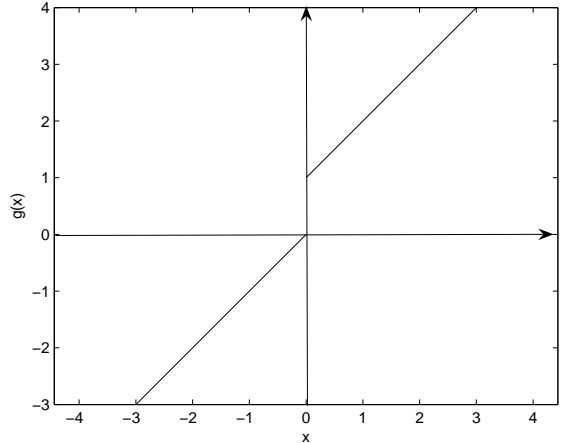}
\caption{An example of function $g\in\mathcal{A}$}
\label{figGraphOfG}
\end{figure}

In the following we first define a special class of discontinuous
functions which will be used throughout this paper.
\begin{definition}(Function class $\mathcal{A}$)\label{ClassAofFunctions}
A function $g$: $\mathbb{R}\mapsto \mathbb{R}$ belongs to
$\mathcal{A}$, denoted by $g\in \mathcal{A}$, if :
\begin{enumerate}
\item $g$ is continuous on $\mathbb{R}$ except for a set with zero
measure, and on each finite interval, the number of discontinuous
points of $g$ is finite.
\item On each interval where $g$ is continuous, $g$ is strictly
increasing;
\item If $x_{0}$ is a discontinuous point of $g$, let $g(x_{0}^{+})=\lim_{x\to
x_{0}^{+}}g(x)$, $g(x_{0}^{-})=\lim_{x\to x_{0}^{-}}g(x)$, then
$g(x_{0}^{+})>g(x_{0}^{-})$.
\end{enumerate}
\end{definition}
\begin{example}\label{examG}
Let
\begin{eqnarray}
g(x)=\left\{
\begin{array}{cc}
x+1,& x>0\\
x,& x<0
\end{array}
\right.,
\end{eqnarray}
then $g\in \mathcal{A}$ with $x=0$ being the only discontinuous point of $g$.
The graph of $g$ is shown in Fig. \ref{figGraphOfG}.
\end{example}
\begin{definition}(shrinking condition)
An absolutely continuous function
$x(t)=[x_{1}(t),\cdots,x_{n}(t)]^{T}$: $\mathbb{R}^{+}\mapsto
\mathbb{R}^{n}$ is {\bf shrinking} if $\max_{i}\{x_{i}(t)\}$ is
nonincreasing and $\min_{i}\{x_{i}(t)\}$ is nondecreasing with
respect to $t$. Furthermore, $x(t)$ is {\bf completely shrinking}£¬
if $x(t)$ is shrinking and
\begin{eqnarray*}
\lim_{t\to+\infty}\max_{i}\{x_{i}(t)\}-\min_{i}\{x_{i}(t)\}=0.
\end{eqnarray*}
\end{definition}
\begin{remark}
It is obvious that if $x(t)$ is shrinking, then the limits of
$\max_{i}\{x_{i}(t)\}$ and $\min_{i}\{x_{i}(t)\}$ exist as $t\to\infty$.
\end{remark}

%Another important
%property for the shrinking condition is given in the following
%
%\begin{lemma}\label{lemPropertyShrinkingCondition}
%Let $\{x^{k}(t)\}$ be a sequence of absolutely continuous functions
%such that $x^{k}(t)$ is (complete) shrinking and $x^{k}(0)=x_{0}$
%for each $k$. If $\{x^{k}(t)\}$ uniformly converge to a function
%$\{x(t)\}$, then $x(t)$ is (completely) shrinking. Moreover, if
%there exists $a\in\mathbb{R}$ such that
%$\lim_{t\to+\infty}x_{i}^{k}(t)=a$ for each $k$, then we have
%\begin{eqnarray}
%\lim_{t\to+\infty}x_{i}(t)=a, \qquad i=1,\cdots,n.
%\end{eqnarray}
%\end{lemma}
%
%{\em Proof. }First, since $x^{k}(t)$ uniformly converge to $x(t)$
%and $x^{k}(t)$ is absolutely continuous, then $x(t)$ is absolutely
%continuous.
%
%If each $x^{k}(t)$ satisfies the complete shrinking condition, then
%$x^{k}(t)$ satisfies the shrinking condition, so from the first
%conclusion, $x(t)$  is shrinking. Then,
%$\lim_{t\to+\infty}\max_{i}{x_{i}^{k}(t)}$ and
%$\lim_{t\to+\infty}\max_{i}{x_{i}(t)}$ exist. So we have
%\begin{align}
%&\lim_{t\to+\infty}\max_{i}\{x_{i}(t)\}
%=\lim_{t\to+\infty}\lim_{k\to+\infty}\max_{i}\{x^{k}_{i}(t)\}
%=\lim_{k\to+\infty}\lim_{t\to +\infty}\max_{i}\{x_{i}^{k}(t)\}\nonumber\\
%=&\lim_{k\to+\infty}\lim_{t\to+\infty}\min_{i}\{x_{i}^{k}(t)\}
%=\lim_{t\to+\infty}\lim_{k\to+\infty}\min_{i}\{x_{i}^{k}(t)\}
%=\lim_{t\to+\infty}\min_{i}\{x_{i}(t)\}
%\end{align}
%and this concludes the proof of the second claim. The remainings are
%obvious. \qquad$\square$

\begin{definition}\label{defSetValuedMap}(Aubin \& Frankowska, 1990)
Let $X$, $Y$ be metric spaces, A map $F$ defined on $E\subseteq X$
is called a set-valued map, if to each $x\in E$, there corresponds a
set $F(x)\subseteq Y$. A set-valued map $F$ is said to be upper
semicontinuous at $x_{0}\in E$ if for any opening set $N$ containing
$F(x_{0})$, there exists a neighborhood $M$ of $x_{0}$ such that
$F(M)\subset N$. $F$ is said to have closed (convex, compact) image,
if for each $x\in E$, $F(x)$ is closed (convex, compact,
respectively).
\end{definition}

\begin{definition}\label{defBasicConditions}(Filippov, 1988)
A set valued map $F$: $\mathbb{R}^{n}\mapsto 2^{\mathbb{R}^{n}}$ is
said to satisfy the basic conditions in a domain $G\subseteq
\mathbb{R}^{n}$ if for any $x\in G$, $F(x)$ is non-empty, bounded,
closed and convex, and $F$ is upper semicontinuous in $x$.
\end{definition}
As to the existence of Filippov solutions, we have the following

\begin{lemma}\label{lemExistenceFilippovSolution}(Filippov,1988)
If a set-valued map $F(x)$ satisfies the basic conditions in the
domain $D\subseteq \mathbb{R}^{n}$, then for any point $x_{0}\in D$,
there exists a solution in $D$ of the following differential
inclusion:
\begin{eqnarray}\label{eqnSolutionDifferentialInclusion}
\dot{x}(t)\in F(x(t)), \qquad x(t_{0})=x_{0}
\end{eqnarray}
over an interval $[t_{0}, t')$ for some $t'>t_{0}$. Moreover, if $F$
satisfies the basic conditions in a closed bounded domain $D$, then
each solution of the differential inclusion
(\ref{eqnSolutionDifferentialInclusion}) lying within $D$ can be
continued either unboundedly as $t$ increases (and decreases), i.e.,
as $t\to \infty$, or until it reaches the boundary of the domain
$D$.
\end{lemma}

\begin{lemma}(Filippov,1988)\label{lemBoundnessOnCompactSet}
Let a set-valued map $F(x)$ be upper semicontinuous on a compactum $K$
and let for each $x\in K$ the set $F(x)$ be bounded, then $F$ is
bounded on $K$.
\end{lemma}
\begin{remark}
It is clear from lemma \ref{lemBoundnessOnCompactSet} that if $F$
satisfies the basic conditions on some compact set $K$, then $F$ is
bounded on $K$.
\end{remark}

\begin{lemma}\label{lemTransformationOnConvexSets}(Filippov,1988)
If $M$ is a bounded closed set and if a function $f$ is continuous,
then the set $f(M)=\{f(x):x\in M\}$ is closed. If $M$ is convex,
$f(x)=Ax+b$, then the set $f(M)=AM+b$ is convex.
\end{lemma}
\begin{remark}
It can be seen from lemma \ref{lemTransformationOnConvexSets} that
if a set-valued map $F(x)$ satisfies the basic condition, then for
any $n\times n$ matrix $T$, the set-valued map $TF(x)=\{Ty:y\in
F(x)\}$ also satisfies the basic condition.
\end{remark}

%\begin{definition}(Filippov,1988)
%Given some $\delta>0$, a vector function $y(t)$ is called a
%$\delta$-solution of a differential inclusion
%\begin{eqnarray}\label{eqnDifferentialInclusionDeltaSolution}
%\dot{x}\in F(x)
%\end{eqnarray}
%\end{definition}
%with $F$ upper semicontinuous in $x$, if $y(t)$ is absolutely
%continuous on a given interval and almost everywhere
%\begin{eqnarray}
%\dot{y}(t)\in F^{\delta}(y(t))
%\end{eqnarray}
%where $y^{\delta}=\{y_{1}\big|\|y_{1}-y\|\le \delta\}$ and
%$F^{\delta}(y)=[\text{co} F(y^{\delta})]^{\delta}$.
%
%\begin{lemma}(Filippov,1988)\label{lemConvergenOfDeltaSolution}
%If a set-valued map $F(x)$ satisfies the basic conditions on an open
%domain $G\subset \mathbb{R}^{n}$, then the limit $x(t)$ of any
%uniformly convergent sequence of $\delta_{k}$-solutions
%$(\delta_{k}\to 0^{+}, k=1,2,\cdots)$ of the inclusion
%(\ref{eqnDifferentialInclusionDeltaSolution}) is a solution of this
%inclusion if the graph of the limiting function $x(t)$ lies within
%$G$.
%\end{lemma}
The following lemma is a generalization of LaSalle invariance principle
for discontinuous differential equations.
\begin{lemma}(Cort{\'e}s, 2006)\label{lemInvariancePrinciple}
Let $V$: $\mathbb{R}^{n}\mapsto \mathbb{R}$ be a locally Lipschitz
and regular function, let $x_{0}\in S\subset \mathbb{R}^{n}$ where
$S$ is compact and strongly invariant with respect to
(\ref{eqnFilippovDefinition}). Assume that either
$\max\tilde{\mathcal{L}}_{f}V(x)\le 0$ or
$\tilde{\mathcal{L}}_{f}V(x)=\emptyset$ for all $x\in S$. Let
$Z_{f,V}=\{x\in\mathbb{R}^{n}|0\in \tilde{\mathcal{L}}_{f}V(x)\}$.
Then, any solution $x(t)$ starting from $x_{0}$ converges to the
largest invariant set $M$ contained in $\overline{Z}_{f,V}\cap S$.
\end{lemma}

\subsection{Probability theory}
Let $\Prb$ denote the probability, and $\Exp$ be the mathematical expectation.
The following are the second Borel-Cantelli Lemma concerning an independent sequence.
\begin{lemma}(Durrett, 2005)\label{lemBorelCantelli}
If the events $\{A_{n}\}$ are independent, then $\sum\Prb\{A_{n}\}=\infty$ implies $\Prb\{A_{n}\text{~~i.o.}\}=1$, where i.o. means infinitely often.
\end{lemma}

\section{Consensus analysis}\hspace{1.25em}\label{secMainResults}
In this section, we will discuss consensus in a network under
nonlinear discontinuous protocols with both fixed topology and switching topologies.

\subsection{Consensus in networks with fixed topology.}

Consider the following consensus protocol in a network of
multiagents with fixed graph topologies:
\begin{eqnarray}\label{sysContinuousFixedTopology}
\dot{x}_{i}=-\sum_{j=1}^{n}l_{ij}g(x_{j}),
\end{eqnarray}
where $g\in\mathcal{A}$ and $L=[l_{ij}]$ is the graph Laplacian.

Denote $\Phi(x)=[\Phi_{1}(x),\cdots,\Phi_{n}(x)]^{\top}$ with
$\Phi_{i}(x)=-\sum\limits_{j=1}^{n}l_{ij}g(x_{j})$, then we can define a
set-valued map
$\mathcal{K}[\Phi_{i}](x)=-\sum_{j=1}^{n}l_{ij}\gamma_{j}$,
with $\gamma_{j}\in \mathcal{K}[g](x_{j})$,
where $\mathcal{K}[g](z)=g(z)$ if $g$ is continuous at $z$,
and $\mathcal{K}[g](z)=[g(z^{-}),g(z^{+})]$ otherwise. Since for any
$x=[x_{1},x_{2},\cdots,x_{n}]^{\top}\in \mathbb{R}^{n}$, the set
$\{\gamma=[\gamma_{1},\gamma_{2},\cdots,\gamma_{n}]^{\top}: \gamma_{i}\in
\mathcal{K}[g](x_{i}),i=1,2,\cdots,n.\}$ is closed and convex, from Lemma
\ref{lemTransformationOnConvexSets}, $\mathcal{K}[\Phi](x)$ is a closed
convex set. The Filippov
solution $x(t)$ to (\ref{sysContinuousFixedTopology}) is defined as
the following differential inclusion:
\begin{eqnarray}\label{eqnDifferentialInclusionDeltaSolution}
\dot{x}_{i}(t)\in\mathcal{K}[\Phi_{i}](x(t)),\qquad
a.e.~~t.
\end{eqnarray}

First, we have the following lemma which says that all the Filippov
solutions of (\ref{sysContinuousFixedTopology}) is shrinking.
\begin{lemma}\label{lemShrinkingOfFilippovSolution}
For any initial value $x_{0}\in\mathbb{R}^{n}$, the Filippov
solution exists and is shrinking, thus, all the solutions can be
extended to $[0,+\infty)$.
\end{lemma}
\noindent{\bf Proof:} It is clear that the set-valued map
$\mathcal{K}[\Phi](x)=\sum_{j=1}^{n}l_{ij}\mathcal{K}[g](x_{j})$
satisfies the basic conditions on any
bounded region of $\mathbb{R}^{n}$, which implies that for any
initial value $x_{0}\in \mathbb{R}^{n}$, the Filippov solution
exists on the interval $[0,t_{1})$ for some $t_{1}>0$ .

Denote $V^{*}(x)=\max_{i}\{x_{i}\}$, $V_{*}(x)=\min_{i}\{x_{i}\}$.
It is easy to see that $V^{*}(x)$ is locally Lipschitz and convex.
In fact, for $x=[x_{1},\cdots,x_{n}]^{\top}$,
$y=[y_{1},\cdots,y_{n}]^{\top}$, and $\lambda\in[0,1]$, we have
\begin{eqnarray*}
|V^{*}(x)-V^{*}(y)|=|\max_{i}\{x_{i}\}-\max_{i}\{y_{i}\}|
\le\max_{i}|x_{i}-y_{i}|
\end{eqnarray*}
and %To see that $V^{*}$ is regular, we only need to verify that
%$V^{*}$ is convex (see \cite{FClarke}). Let $\lambda\in[0,1]$, we
%have
\begin{eqnarray*}
V^{*}(\lambda x+(1-\lambda)y)&=&\max_{i}\{\lambda
x_{i}+(1-\lambda)y_{i}\}\\
&\le&\lambda\max_{i}\{x_{i}\}+(1-\lambda)\max_{i}\{y_{i}\}\\
&=&\lambda V^{*}(x)+(1-\lambda)V^{*}(y),
\end{eqnarray*}
Therefore, $V^{*}$ is regular and
\begin{eqnarray*}
\frac{dV^{*}(x(t))}{dt}\in \tilde{\mathcal{L}}_{\Phi}V^{*}(x),\qquad
a.e. ~t.
\end{eqnarray*}
where $\tilde{\mathcal{L}}_{\Phi}V^{*}$ is the set-valued Lie
derivative of $V^{*}$ with respect to $\Phi$.

We will prove that $V^{*}(x(t))$ is nonincreasing and $V_{*}(x(t))$
is nondecreasing. Here, we only show that $V^{*}(x(t))$ is
nonincreasing, and a similar argument can apply to $V_{*}(x(t))$.

Now, we will prove that for each $x$, either
$\tilde{\mathcal{L}}_{\Phi}V^{*}(x)=\emptyset$ or
$\max\{\tilde{\mathcal{L}}_{\Phi}V^{*}(x)\}\le 0$. Given
$x=[x_{1},\cdots,x_{n}]^{\top}\in\mathbb{R}^{n}$, let
%$\overline{I}_{x}$ be the set of indices defined as
$\overline{I}_{x}=\{i\in\{1,\cdots,n\}: x_{i}=\max_{j}\{x_{j}\}\}$.
We have $\partial V^{*}(x)=\rm{co}\{e_{i}: i\in \overline{I}_{x}\}$. If $a\in
\tilde{\mathcal{L}}_{\Phi}V^{*}(x)$, then there exists some
$v=[v_{1},\cdots,v_{n}]^{\top}\in \mathcal{K}[\Phi](x)$ such that
$a=v\cdot \zeta$ for each $\zeta\in\partial V^{*}(x)$. Therefore,
$v_{i}=a$ for $i\in\overline{I}_{x}$.

Noting
\begin{eqnarray*}
v_{i}= -\sum_{j=1}^{n}l_{ij}\gamma_{j}=-\sum_{j=1,\atop{j\ne i}}^{n}l_{ij}(\gamma_{j}-\gamma_{i}),
\end{eqnarray*}
for some $\gamma_{j}\in\mathcal{K}[g](x_{j})$, if $g$
is continuous at $x_{i}$, then $\gamma_{j}=g(x_{j})=g(x_{i})=\gamma_{i}$ for $j\in \overline{I}_{x}$,
and $\gamma_{j}<\gamma_{i}$ for $j\not\in\overline{I}_{x}$. So in this case we have $v_{i}\le 0$. Otherwise,
$g$ is discontinuous at $x_{i}$. If $a>0$, then for each $i\in \overline{I}_{x}$, $v_{i}=a>0$.
Let $\overline{i}\in \overline{I}_{x}$ be one index satisfying $\gamma_{\overline{i}}=\max\{\gamma_{i}:~i\in \overline{I}_{x}\}$. Then we obviously have $v_{\overline{i}}\le 0$, which is a contradiction. So in this case
we also have $a\le 0$.

From Lemma \ref{lemEvolutionAlongFilippovSolution},
\begin{equation*} \frac{dV^{*}(x(t))}{dt}\le 0,\quad ~~\text{a.e.}~~t.\end{equation*}
Thus $V^{*}(x(t))$ is nonincreasing. A similar argument can show that
$V_{*}(x(t))$ is nondecreasing. So $x(t)$ is shrinking.
The second claim then directly follows from Lemma \ref{lemExistenceFilippovSolution}.
 \qquad$\square$

Based on lemma \ref{lemShrinkingOfFilippovSolution}, we can prove
following theorem concerning the consensus of system
(\ref{sysContinuousFixedTopology}).

\begin{theorem}\label{thmMainFixTopology}
The system (\ref{sysContinuousFixedTopology}) will achieve consensus
for any initial value if and only if the graph of $L$ has spanning
trees. And the consensus value is $\Wra(x(0),L)$. Furthermore, if the graph of $L$
is strongly connected, and $g$ is discontinuous at $\Wra(x(0),L)$,
then finite time convergence can be achieved.
\end{theorem}
{\bf \indent Proof:} See Appendix \ref{appenProofThm1}.

It can be seen that Theorem \ref{thmMainFixTopology} is quite
similar to the result obtained in literature for continuous
consensus protocols. So the protocol
\eqref{sysContinuousFixedTopology} can be seen as natural extensions
of the continuous protocols. Intuitively, if a networks has spanning
trees, then the information from the roots can be sent to all other
nodes in the network. And the roots can exchange information with
each other. So the network can finally reach a consensus. If a
network has no spanning trees, from the proof of Theorem
\ref{thmMainFixTopology}, there are two possible cases. Case I:
there exists an isolated subgraph that has no connection with other
parts of the network. In this case the isolated subgraph can not
exchange information with other parts of the network, and consensus
can not be reached. Case II: there are no isolated subgraphs. In
this case, the network has a subgraph that has spanning trees. There
are edges from nodes outside this subgraph to nodes of this subgraph
which are not roots. Fig. \ref{figNoSpanGraph} provides an example.
In this case, the roots in the subgraph can not exchange information
with nodes outside the subgraph, since they can neither send their
information to the nodes outside the subgraph, nor receive
information from nodes outside the subgraph. As a result, consensus
also can not be reached. In the following, we will provide some
examples to illustrate the theoretical results.

\begin{example}
The graph shown in Fig. \ref{figSpanGraph} may be called a {\bf ``double-star" graph}. It has spanning trees, with $\{v_{1},v_{2}\}$ being the set of roots. Yet this graph is not strongly connected. If we take the weight of each edge to be $1$, then the graph Laplacian is
\begin{eqnarray*}
L=\left[
\begin{array}{cc}
L_{11} & 0\\
L_{21} &I_{10}
\end{array}
\right],\\[1pt]
\end{eqnarray*}
with $L_{11}=\left[\begin{array}{cc}1&-1\\-1&1\end{array}\right]$, $L_{21}=\left[\begin{array}{cccccccccc}
-1&-1&-1&-1&-1&0&0&0&0&0\\
0&0&0&0&0&-1&-1&-1&-1&-1
\end{array}\right]^{\top}$, and $I_{10}$ being the $10\times 10$ identity matrix.
For any $x=[x_{1},x_{2},\cdots,x_{12}]^{\top}\in \mathbb{R}^{12}$, $\Wra(x,L)=(x_{1}+x_{2})/2$. The simulation result is provided in Fig.
\ref{figFix}, where $g$ is given in Example \ref{examG}, and the initial value $x(0)$ is randomly chosen. The
position of $\Wra(x(0),L)=(x_{1}(0)+x_{2}(0))/2$ is labeled
on the right side with a `+'. It can be seen that the agents finally
reach a consensus on $\Wra(x(0),L)$.

\begin{figure}
\centering
\includegraphics[width=0.5\textwidth]{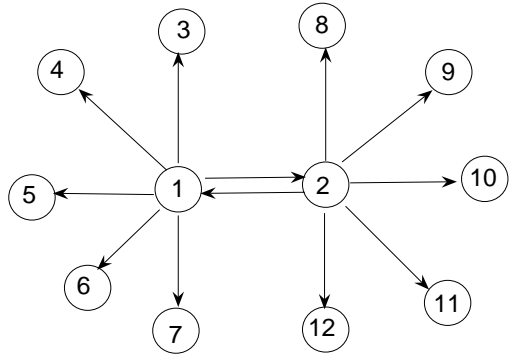}
\caption{A ``double star" network}
\label{figSpanGraph}
\end{figure}
\begin{figure}[!t]
\begin{center}
\includegraphics[width=0.5\textwidth]{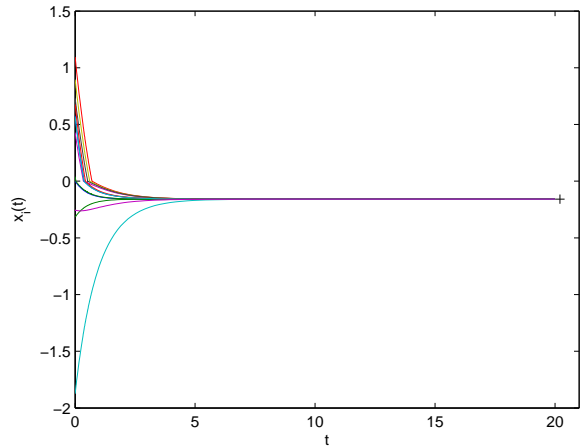}
\caption{Consensus in a ``double-star" network}
\label{figFix}
\end{center}
\end{figure}
\end{example}

\begin{example}
Fig. \ref{figNoSpanGraph} provides an example of a graph that has no spanning trees. This graph has no isolated subgraphs.
The subgraph induced by
$\{v_{1},v_{2},v_{3},v_{6}\}$ has spanning trees, with $\{v_{1},v_{2}\}$ being the root set. And there are
edges from $\{v_{4},v_{5}\}$ to $\{v_{3},v_{6}\}$. So this graph belongs to the second case discussed above.
And it can not reach a consensus for arbitrary initial value. For each edge, we take the weight as $1$. Then
the graph Laplacian is
\begin{eqnarray*}
L=\left[
\begin{array}{rrrrrr}
1 & -1 & 0 & 0 & 0 & 0 \\
-1& 1  & 0 & 0 & 0 & 0 \\
-1 &0  & 2 & 0 &-1 & 0 \\
0 &-1  & 0 & 2 & 0 &-1 \\
0 & 0  & 0 & 0 & 1 &-1 \\
0 & 0  & 0 & 0 & -1& 1
\end{array}
\right]
\end{eqnarray*}
The simulation results are presented in Fig. \ref{figNonConsensus}, with $g$ being given in Example \ref{examG}. It can be seen that no consensus is realized.
\begin{figure}
\centering
\includegraphics[width=0.5\textwidth]{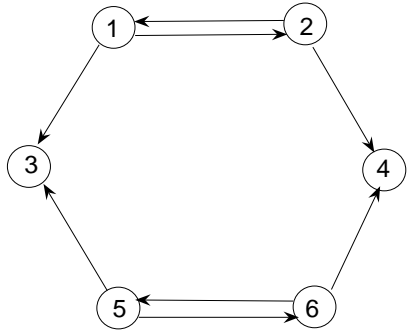}
\caption{An example of a graph that has no spanning trees}
\label{figNoSpanGraph}
\end{figure}
\begin{figure}
\centering
\includegraphics[width=0.5\textwidth]{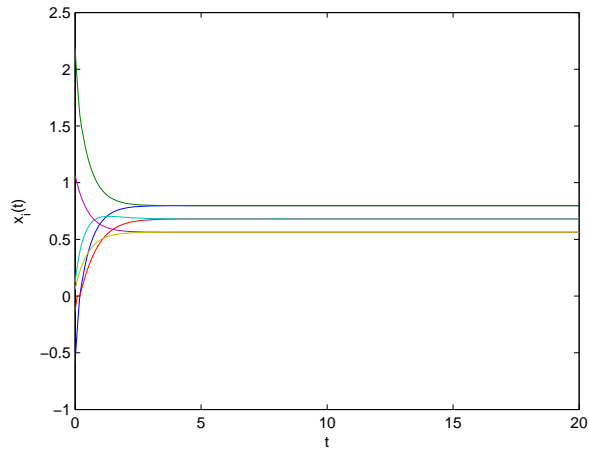}
\caption{No consensus can be reached in the graph of Fig. \ref{figNoSpanGraph}.}
\label{figNonConsensus}
\end{figure}
\end{example}

\subsection{Consensus in networks with randomly switching topologies}

In this section, we will investigate consensus in networks of
multiagents under nonlinear protocols over graphes with randomly switching
topologies.

Consider the following dynamical system:
\begin{eqnarray}\label{sysContinuousSwitchingTopology}
\dot{x}_{i}(t)=-\sum_{j=1}^{n}l^{k}_{ij}g(x_{j})&&t\in
[t_{k},t_{k+1}),
\end{eqnarray}
where $g\in\mathcal{A}$ and $L^{k}=[l_{ij}^{k}]$ is the graph
Laplacian for the underlying graph on the time interval
$[t_{k},t_{k+1})$. At each time point $t_{k}$ there is a switching
of the network topology. We consider the case that $L^{k}$ is a random sequence.
%Every graph Laplacian $L^{k}=[l^{k}_{ij}]$
%is generated independently in the following way: Given a weight
%matrix $W=[w_{ij}]$, and some $p\in (0,1)$, for $i\ne j$
%\begin{eqnarray*}
%l^{k}_{ij}=\left\{
%\begin{array}{cc}
%w_{ij}&\textnormal{ with probability } p\\
%0& \textnormal{ with probability } 1-p.
%\end{array}
%\right.
%\end{eqnarray*}
Denote $\Delta t_{k}=t_{k+1}-t_{k}$, in this following, we make

\begin{assumption}\label{asumpIndependentConditions}
\begin{enumerate}
\item $\{\Delta t_{j}\}$ is independent and identically distributed;
%such that $E\Delta t_{k}<+\infty$ and $E\Delta t_{k}^{2}<+\infty$.
%\item $\{L^{k}\}$ is independent and identically distributed.
\item the sequence $\{\Delta t_{i}\}$ and $\{L^{k}\}$ are independent;
\item $\{L^{k}\}$ is uniformly bounded.
\end{enumerate}
\end{assumption}

\begin{assumption}\label{asumpSeperateConditions}
There exists $\varepsilon>0$ such that for any
$\alpha,\beta\in\mathbb{R}$ with $\alpha\ne\beta$ and $g$ is
continuous at $\alpha$, $\beta$, it satisfies that
\begin{equation*}
\frac{g(\alpha)-g(\beta)}{\alpha-\beta}\ge \varepsilon.
\end{equation*}
\end{assumption}
\begin{remark}
It is easy to verify that under Assumption
\ref{asumpSeperateConditions}, for any $\alpha,\beta\in\mathbb{R}$
with $\alpha\ne \beta$ and $v_{1}\in \mathcal{K}[g](\alpha)$,
$v_{2}\in\mathcal{K}[g](\beta)$, it satisfies that
\begin{eqnarray*}
\frac{v_{1}-v_{2}}{\alpha-\beta}\ge \epsilon.
\end{eqnarray*}
\end{remark}

First, we will prove the following Theorem for almost sure consensus.
\begin{theorem}\label{thmMain}
Under Assumption \ref{asumpIndependentConditions},
\ref{asumpSeperateConditions}, the system
(\ref{sysContinuousSwitchingTopology}) will achieve consensus almost
surely if there exists $\delta>0$ such that
\begin{eqnarray*}
\Prb\{\mathcal{G}(L^{k})\text{ is $\delta$-scrambling for infinitely }k\}=1.
\end{eqnarray*}
\end{theorem}
\noindent{\bf Proof:} See Appendix \ref{appenProofThm2}.

From Theorem \ref{thmMain}, we can have the following corollary concerning switching sequence
$\{L^{k}\}$ which is independent and identically distributed.
\begin{corollary}\label{corIID}
Under Assumption \ref{asumpIndependentConditions},\ref{asumpSeperateConditions}, if $\{L^{k}\}$ is independent and
identically distributed, then the system
(\ref{sysContinuousSwitchingTopology}) will achieve consensus almost
surely if $\Exp \eta(-L^{k})>0$.
\end{corollary}
\noindent{\bf Proof:} Denote $\delta=\Exp \eta(-L^{k})>0$, and $M=\sup \eta(-L^{k})<\infty$. Then
\begin{eqnarray*}
\delta=\Exp \eta(-L^{k})\le \frac{\delta}{2}\Prb\{\eta(-L^{k})\le  \frac{\delta}{2}\}+M\Prb\{\eta(-L^{k})>\frac{\delta}{2}\}\le \frac{\delta}{2}+M\Prb\{\eta(-L^{k})>\frac{\delta}{2}\},
\end{eqnarray*}
which implies
\begin{eqnarray*}
\Prb\{\eta(-L^{k}>\delta/2)\}\ge \frac{\delta}{2M}.
\end{eqnarray*}
From the second Borel-Cantelli Lemma (Lemma \ref{lemBorelCantelli}), we have
\begin{eqnarray*}
\Prb\{\eta(-L^{k})>\frac{\delta}{2} \text{ infinitely often}\}=1.
\end{eqnarray*}
The conclusion follows from Theorem \ref{thmMain}.

\section{Applications to a generalized blinking model}\label{secNumericalSimulations}
In this section, we will show how the theoretical results can be applied to analyze real world network models. For this purpose, we consider a generalized blinking network model.

The original blinking model was proposed in Belykh, Belykh, \& Hasler (2004). It is a kind of small world networks
that consists of a regular lattice of cells with constant $2K$ nearest neighbor couplings and time dependent
on-off couplings between any other pair of cells. In each time interval of duration $\tau$ each time dependent
coupling is switched on with a probability $p$, and the corresponding switching random variables are independent
for different links and for different times. It is a good model for many real-world dynamical networks such as
computers networked over the Internet interact by sending packets of information, and neurons in our brain
interact by sending short pulses called spikes, etc.

On the other hand, this model is still quite restrictive in several aspects. First, this model is an undirected model. Second, the
duration between any two successive switchings may not be identical, nor may it be small sometimes. And it may even be not deterministic, but just a random variable. Finally, the basic regular $2K$ nearest neighbor coupling lattice may not exist, or we can say $K=0$ in such case.

Based on the above analysis, we make the following generalizations on the original blinking model.
First, we assume the model to be a directed graph. For every two vertices $v_{i}$, $v_{j}$ that have random switching links between them, the switching of the edge from $v_{i}$ to $v_{j}$ is independent of that from $v_{j}$ to $v_{i}$. Second, we assume the duration between every two successive switchings is a random variable, and each duration is independent of others. Finally, we assume that $K$ may be zero in the basic $2K$ nearest neighbor lattice. That is, no links exist with probability $1$.

It is obvious that in this generalized model, the sequence of the durations are independent and identically distributed. And the underlying graph sequence $\{\mathcal{G}^{k}\}$ is also independent and identically distributed. For each $\mathcal{G}^{k}$, since different links are switched on independently, it is obvious that there is a positive probability that $\mathcal{G}^{k}$ is a complete graph. Since a complete graph is scrambling, if we set the weight of each link to be $\delta>0$, then $\mathcal{G}^{k}$ is $\delta$-scrambling for some $\delta>0$ with a positive probability. From Corollary \ref{corIID}, we can see that the discontinuous consensus protocol \eqref{sysContinuousSwitchingTopology} will realize consensus almost surely on a generalized blinking model.

\begin{figure}
\begin{center}
\includegraphics[width=0.5\textwidth]{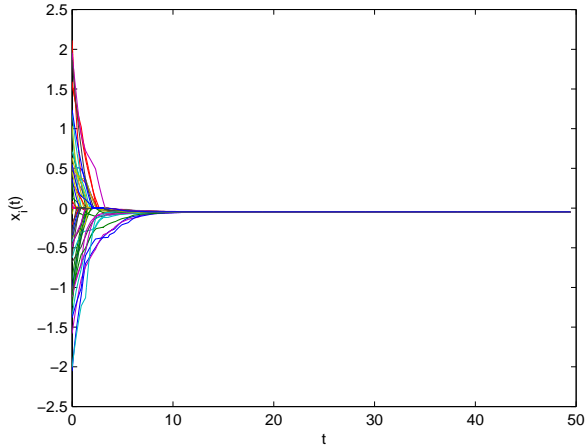}%[width=2.25in]
\caption{Consensus in a generalized blinking model.} \label{figBlinking}
\end{center}
\end{figure}

In the simulation, we choose a network with $50$ nodes, $K=0$, $p=0.1$, and the weight of each link to be $0.1$. The duration between every two successive switching is a random variable uniformly distributed on $(0,1)$. Let $g$ be as given in Example \ref{examG}. The initial value is chosen randomly. The simulation results are presented in Fig. \ref{figBlinking}. It can be seen that consensus can be reached almost surely.

\section{Conclusion}\hspace{1.25em}\label{secConclusions}
In this paper, we investigate consensus in networks of
multiagents under nonlinear discontinuous protocols. First, we consider networks with fixed topology described by weighted directed graphs. Compared to existing results concerning discontinuous consensus protocols, we do not require the underlying graph to be strongly connected.
Instead, we prove that a directed spanning tree is sufficient and necessary to realize consensus.
And we can also locate the consensus value. This result can be seen as an extension of continuous protocols if we take continuous protocols as special case of discontinuous ones. Under this viewpoint, we establish a more generalized theoretical framework for consensus analysis. Second, we consider networks with randomly switching topologies. We provide sufficient conditions for the network to achieve consensus almost surely based on the scramblingness of the underlying graphs.
Particularly, we consider the case when the switching sequence is independent and identically distributed. Compared to existing results on discontinuous protocols, we do not require the network to be connected at each time point. Finally, as application of the theoretical results, we study a generalized blinking model and show that consensus can be realized almost surely under the proposed discontinuous protocols.

%\begin{thebibliography}{99}
\noindent{\bf\Large References}
\begin{description}

\item %1
Aubin, J., Frankowska, H. (1990). Set-valued analysis, Boston: Birkhauser.

\item
Belykh, I., Belykh, V., \& Hasler, M. (2004). Blinking model and synchronization in small world networks
with a time-varying coupling, Physica D, 195, 188-206.

%\item %24
%Chen, T., Liu, X., \& Lu, W. (2007). Pinning complex networks by a
%single controller, {\it IEEE Trans. Circuits Syst.-I}, 54,
%1317-1326.
%
%\item %23
%Chen, T., \& Zhu, Z. (2007). Exponential synchronization of
%nonlinear coupled dynamical networks, {\it Int. J. Bifur. Chaos},
%17, 999-1005.

\item %26
Clarke, F. (1983). Optimization and nonsmooth analysis, New York:
Wiley.

\item %2
Cort{\'e}s, J. (2006). Finite-time convergent gradient flows with
applications to network consensus, {\it Automatica}, 42, 1993-2000.

\item %6
Cort{\'e}s, J. (2008). Distributed algorithms for reaching consensus
on general functions, {\it Automatica}, 44, 726-737.

\item %1
Cort{\'e}s, J., \& Bullo, F. (2005). Coordination and geometric
optimization via distributed dynamical systems, {\it SIAM Journal on Control
and Optimization}, 44, 1543-1574.

\item %27
Durrett, R. (2005). Probability: Theory and Examples, 3rd ed.
Belmont, CA: Duxbury Press.

\item %4
Fax, A., \& Murray, R. (2004). Information flow and cooperative
control of hehicle formations, {\it IEEE Transactions on Automatic Control},
49, 1465-1476.

\item %25
Filippov, A. (1988). Differential equations with discontinuous
righthand sides, Kluwer.

\item %5
Hui, Q., Haddad, W., \& Bhat, S. (2008). Semistability theory for
differential inclusions with applications to consensus problems in
dynamical networks with switching topology, {\it 2008 American
Control Conference}, 3981-3986.

\item %13
Liu, B., \& Chen, T. (2008). Consensus in networks of multiagents
with cooperation and competetion via stochastically switching
topoloiges, {\it IEEE Transactions on Neural Networks}, 19, 1967-1973.

\item
Liu, B., Lu, W., \& Chen, T. (2011). Consensus in networks of multiagents with
switching topologies modeled as adapted stochastic processes, {\it SIAM Journal
on Control and Optimization}, 49, 227-253.

\item %14
Liu, X., Chen, T., \& Lu, W. (2009). Consensus problem in directed
networks of multi-agents via nonlinear protocols, {\it Physics
Letters A}, 373, 3122-3127.

%\item %22
%Lu, W., Atay, F., \& Jost, J. (2007). Synchronization of
%discrete-time dynamical networks with time-varying couplings, {\it
%SIAM J. Math. Anal.}, 39, 1231-1259.

\item %15
Lu, W., \& Chen, T. (2004). Synchronization of coupled connected
neural networks with delays, {\it IEEE Transactions on Circuits and Systems-I}, 51,
2491-2503.

\item
Matei, I., Martins, N., \& Baras, J. (2008). Almost sure
convergence to consensus in Markovian random graphs, {\it Proceedings of
the 47th IEEE Conference on Decision and Control}, Cancun, Mexico, 3535-3540.

\item %7
Moreau, L. (2004). Stability of continuous-time distributed
consensus algorithms, {\it 43rd IEEE Conference on Decision and
Control}, 3998-4003.

%\item %8
%Moreau, L. (2006). Stability of multiagent systems with
%time-dependent communication links, {\it IEEE Trans. Automat.
%Control}, 51, 169-182.

%\item %19
%Nishikawa, T., \& Motter,  A. (2006). Synchronization is optimal in
%nondiagonalizable networks, {\it Phys. Rev. E}, 73, 065106.
%
%\item %20
%Nishikawa, T., \& Motter, A. (2006). Maximum performance at minimum
%cost in network synchronization, {\it Physica D}, 224, 77-89.

\item %3
Olfati-Saber, R. (2006). Flocking for multi-agent dynamical systems:
Algorithms and theory, {\it IEEE Transactions Automatic Control},  51,
401-420.

\item %10
Olfati-Saber, R., Fax, J., \& Murray, R. (2007). Consensus and
cooperation in networked multi-agent systems, {\it Proceedings of IEEE}, 95,
215-233.

\item %9
Olfati-Saber, R., \& Murray, R. (2004). Consensus problems in
networks of agents with switching topology and time delays, {\it
IEEE Transactions on Automatic Control}, 49, 1520-1533.

\item %11
Ren, W., Beard, R., \& Atkins, E. (2005). A survey of consensus
problems in multi-agent coordination, {\it Proceedings of the
American Control Conference, Holland, OR}, 1859-1864.

\item
Salehi, A. \& Jadbabaie, A. (2007). Necessary and sufficient
conditions for consensus over random independent and identically
distributed switching graphs, {\it Proceedings of the 46th IEEE
Conference on Decision and Control}, 4209-4214.

%\item %12
%Ren, W., Beard, R., \& Kingston, D. (2005). Multi-agent Kalman
%consensus with relative uncertainty, {\it Proceedings of the
%American Control Conference, Holland, OR}, 1865-1870.

%\item %16
%Restrepo, J., Hunt, B., \& Ott, E.,  (2005). Onset of
%synchronization in large networks of coupled oscillators, {\it Phys.
%Rev. E}, 71, 036151.
%
%\item %17
%Restrepo, J., Ott, E., \& Hunt, B. (2006). Synchronization in large
%directed networks of coupled phase oscillators, {\it Chaos}, 16,
%015107.

%\item %21
%Sorrentino, F., di Bernardo, M., Huerta Cu\'{e}llar, G., \&
%Boccaletti, S. (2006). Synchronization in weighted scale-free
%networks with degree¨Cdegree correlation, {\it Physica D}, 224,
%123-129.

%\item %18
%Zhou, C., Motter, A., \& Kurths, J. (2006). Universality in the
%Synchronization of Weighted Random Networks, {\it Phys. Rev. Lett.},
%96, 034101.

\end{description}
\appendix
\section{Proof of Theorem \ref{thmMainFixTopology}}\label{appenProofThm1}
{\em Sufficiency: }Let $V=V^{*}-V_{*}$, where
$V^{*}$ and $V_{*}$ are defined as in lemma
\ref{lemShrinkingOfFilippovSolution}. Then $V$ is locally Lipschitz
and regular.

Given any initial value $x(0)\in\mathbb{R}^{n}$, denote
$\overline{x}_{0}=\max_{i}\{x_{i}(0)\}$,
$\underline{x}_{0}=\min_{i}\{x_{i}(0)\}$, and
$S=\{x=[x_{1},\cdots,x_{n}]^{\top}\in\mathbb{R}^{n}:
\underline{x}_{0}\le x_{i}\le\overline{x}_{0}\}$. By lemma
\ref{lemShrinkingOfFilippovSolution}, $S$ is strongly invariant. Let
 $\overline{Z}_{\Phi,V}=\{x\in \mathbb{R}^{n}:0\in
\tilde{\mathcal{L}}_{\Phi}V\}$ and $M$ be the largest weakly
invariant set contained in $\overline{Z}_{\Phi,V}\cap S$. By Lasalle
Invariance Principle (see Lemma \ref{lemInvariancePrinciple}), we have
\begin{eqnarray*}
\Omega(x(t))\subseteq M,
\end{eqnarray*}
where $\Omega(x(t))$ is the positive limit set of $x(t)$.

Let
$\mathscr{C}=\{x=[x_{1}\cdots,x_{n}]^{\top}\in\mathbb{R}^{n},x_{1}=\cdots=x_{n}\}$
be the consensus manifold, we claim that
  ($M\subseteq
\mathscr{C}\cap S$.) $Z_{\Phi,V}\subset \mathscr{C}$. Otherwise,
there exists $x=[x_{1},\cdots,x_{n}]^{\top}\in M$ such that
$\max_{i}\{x_{i}\}>\min_{i}\{x_{i}\}$ and
\begin{eqnarray*}
0=\dot{x}_{i}&\in&\sum_{j=1}^{n}l_{ij}\mathcal{K}[g](x_{j}),
\end{eqnarray*}
which means that for some $v=[v_{1},\cdots,v_{n}]^{\top}$ with
$v_{i}\in \mathcal{K}[g](x_{i})$  and
$\max_{i}\{v_{i}\}>\min_{i}\{v_{i}\}$, and
\begin{eqnarray*}
\sum_{j\in N_{i}}l_{ij}(v_{j}-v_{i})=0,
\end{eqnarray*}

Let $\overline{I}_{v}=\{i:v_{i}=\max_{j}\{v_{j}\}\}$, and
$\underline{I}_{v}=\{i:v_{i}=\min_{j}\{v_{j}\}\}$. First, from the
monotonicity of $g$, we have that $\overline{I}_{v}\subseteq
\overline{I}_{x}$ and $\underline{I}_{v}\subseteq
\underline{I}_{x}$. For $i\in \overline{I}_{v}$, we have
\begin{eqnarray*}
0=\sum_{j\in N_{i}}l_{ij}(v_{j}-v_{i}).
\end{eqnarray*}
This implies that $N_{i}\subseteq \overline{I}_{v}$ for all $i\in
\overline{I}_{v}$. By induction arguments it can be seen that the
root set of the spanning trees is contained in $\overline{I}_{v}$. A
similar argument reveals that the root set of spanning trees are
contained in $\underline{I}_{v}$. But from the assumption that
$\max_{i}\{x_{i}\}>\min_{i}\{x_{i}\}$, we can obviously have that
$\overline{I}_{v}\cap \underline{I}_{v}=\emptyset$, which is a
contradiction.

Based on previous derivation, we proved that $\Omega(x(t))\subseteq
M\subseteq \mathscr{C}$. Next, we will show that $\Omega(x(t))$ only
contains one point. Otherwise, there exist $u=[a,\cdots,a]^{\top}\in
\Omega(x(t))$, $v=[b,\cdots,b]^{\top}\in \Omega(x(t))$, $a\ne b$.
Assume $a>b$. Then there exists a sequence $t_{n}\to+\infty$ as
$n\to+\infty$ such that
$\lim_{n\to+\infty}\min_{i}\{x_{i}(t_{n})\}=a$. By the fact that
$\min_{i}\{x_{i}(t)\}$ is nondecreasing, we have
$\lim_{t\to+\infty}\min_{i}\{x_{i}(t)\}=a$, which implies $b>a$. A
contradiction.

Summing up, we have proved that $\lim_{t\to+\infty}x(t)=x_{\infty}$
for some $x_{\infty}\in \mathscr{C}\cap S$. This completes the proof
of the sufficiency.

{\em Necessity:} Let $\mathcal{G}(L)=\{\mathcal{V},\mathcal{E}\}$ be
the graph of $L$, if $\mathcal{G}(L)$ doesn't have a spanning tree,
then there is a subgraph $\mathcal{G}_{s}$ of $\mathcal{G}(L)$ that
is a maximum spanning tree, i.e., if there exists a subgraph
$\mathcal{G}_{s'}$ of $\mathcal{G}(L)$ that is a spanning tree and
contains $\mathcal{G}_{s}$, then $\mathcal{G}_{s}=\mathcal{G}_{s'}$.
Let $\mathcal{V}_{s}$ be the vertex set of $\mathcal{G}_{s}$, and
$\mathcal{V}_{c}=\mathcal{V}\backslash \mathcal{V}_{s}$. Then,
$\mathcal{V}_{c}\ne \emptyset$. Let $\mathcal{V}_{sr}$ be the set of
roots of $\mathcal{G}_{s}$, and let
$\mathcal{V}_{s'}=\mathcal{V}_{s}\backslash \mathcal{V}_{sr}$.
Obviously, the following properties hold:
\begin{enumerate}
\item There are no edges from $\mathcal{V}_{s}$ to
$\mathcal{V}_{c}$.
\item There are no edges from $\mathcal{V}\backslash
\mathcal{V}_{sr}$ to $\mathcal{V}_{sr}$.
\end{enumerate}
Here, for two vertex sets, an edge from one to the other means an
edge from some vertex in the former to some vertex in the latter.

Then there are two cases to be considered. For simplicity, we denote
each vertex by index, and the vertex should be renumbered if
necessary.
\begin{enumerate}
\item $\mathcal{V}_{s'}=\emptyset$.

In this case, after proper renumbering, from the above mentioned two
properties, the matrix $L$ has the following form:
$L=\left[\begin{array}{cc}L_{1}&0\\0&L_{2}\end{array}\right]$, where
$L_{1}$, $L_{2}$ correspond to $\mathcal{V}_{sr}$,
$\mathcal{V}_{c}$, respectively. Let $n_{1}$ be the dimension of
$L_{1}$ and
$x_{0}=[\underbrace{a,\cdots,a}_{n_{1}},\underbrace{b,\cdots,b}_{n-n_{1}}]^{\top}$
with $a\ne b$, then obviously, $x(t)\equiv x_{0}$ is a solution
which can not achieve any consensus.

\item $\mathcal{V}_{s'}\ne\emptyset$.

In this case, after proper renumbering, from the above mentioned two
properties, the matrix $L$ has the following form:
$L=\left[\begin{array}{ccc}L_{1}&0&0\\ *
&L_{2}&*\\0&0&L_{3}\end{array}\right]$, where $L_{1}$, $L_{2}$,
$L_{3}$ correspond to $\mathcal{V}_{sr}$, $\mathcal{V}_{s'}$,
$\mathcal{V}_{c}$, respectively, and ``*" can be anything. Let
$n_{i}$ be the dimensions of $L_{i}$ for $i=1,2,3$. Let
$x_{0}=[\underbrace{a,\cdots,a}_{n_{1}},b_{1},\cdots,b_{n_{2}},\underbrace{c,\cdots,c}_{n_{3}}]^{T}$
for some $a\ne c$ and $b_{i}\in\mathbb{R}$, then we have
$\dot{x}_{i}\equiv 0$ for $i\in \mathcal{V}_{sr}\cup
\mathcal{V}_{c}$. Therefore, for any solution $x(t)$ starting from
$x_{0}$, it holds that
\begin{eqnarray*}
x_{i}(t)\equiv \left\{\begin{array}{rl}a & i\in \mathcal{V}_{sr}\\
b& i\in \mathcal{V}_{c}.\end{array}\right.
\end{eqnarray*}
no consensus will be achieved.
\end{enumerate}

At last, we prove the consensus value is $\Wra(x(0),L)$.
Suppose that $\mathcal{G}(L)$ has spanning
trees, and $L$ is of the following form
\begin{eqnarray}\label{NormFormofL}
L=\left[
\begin{array}{cc}
L_{1}&0\\
* & L_{2},
\end{array}
\right]
\end{eqnarray}
where $L_{1}$ corresponds to the vertex set of all the roots of the
spanning trees.  In such case, we have that $L_{1}=L_{I_{r}}$.

Let $\xi=[\xi_{1},\cdots,\xi_{\#I_{r}}]^{\top}$ be the eigenvector
corresponding to the zero eigenvalue of $L_{1}$. Assume
$\sum_{i=1}^{\#I_{r}}\xi_{i}=1$, and let
$\overline{x}(t)=\xi^{\top}x_{I_{r}}(t)$, then for almost all $t$,
\begin{eqnarray*}
\dot{\overline{x}}=\sum_{i=1}^{\#I_{r}}\xi_{i}\sum_{j=1}^{\#I_{r}}l_{ij}\gamma_{j}(t)
=\sum_{j=1}^{\#I_{r}}(\sum_{i=1}^{\#I_{r}}\xi_{i}l_{ij})\gamma_{j}(t)
=0,
\end{eqnarray*}
where $\gamma_{j}(t)\in \mathcal{K}[g](x_{j}(t))$.
This implies that $\overline{x}(t)\equiv \overline{x}(0)$. Since
$\lim_{t\to+\infty}x_{i}(t)=x_{\infty}$, we have
$\lim_{t\to+\infty}\sum_{i=1}^{\#I_{r}}\xi_{i}x_{i}(t)=x_{\infty}$.
thus  $x_{\infty}=\overline{x}(0)$.

At last, we prove finite time convergence when $g$ is discontinuous at $\Wra(x(0),L)$.
Denote $\bar{x}=\Wra(x(0),L)$, and let $\bar{\gamma}\in \mathcal{K}[g](\bar{x})$.
For $x=[x_{1},\cdots, x_{n}]^{\top}$, define a function
\begin{eqnarray*}
V_{L}(x)=\sum_{i=1}^{n}\xi_{i}\int_{\bar{x}}^{x_{i}}[g(s)-\bar{\gamma}]ds,
\end{eqnarray*}
where $\xi=[\xi_{1},\cdots,\xi_{n}]^{\top}$ is the positive left eigenvector corresponding
to the zero eigenvalue of $L$ such that $\sum_{i=1}^{n}\xi_{i}=1$. Then it is obvious that
$V_{L}\ge 0$, $V_{L}(x)=0$ if and only if $x_{i}=\bar{x}$ for each $i$. Furthermore, since
$g(x)$ is strictly increasing, and $\xi_{i}>0$, $i=1,\cdots,n$, $V_{L}(x)$ is convex, thus regular.
Also, $V_{L}(x)$ is locally Lipschitz. So from Lemma \ref{lemEvolutionAlongFilippovSolution}, $\frac{dV_{L}(x(t))}{dt}$ exists for a.e. $t$, and
\begin{eqnarray*}
\frac{dV_{L}(x(t))}{dt}=\tilde{\mathcal{L}}_{\Phi}(x(t)),~~\text{a.e.}~~ t.
\end{eqnarray*}
Since from definition, $\partial V_{L}(x)=\{[\gamma_{1},\cdots,\gamma_{n}]^{\top}: \gamma_{i}\in \mathcal{K}[g](x_{i}),~~i=1,\cdots,n\}$, if $\tilde{\mathcal{L}}_{\Phi}(x)\ne \emptyset$, then either
$g$ is continuous at each $x_{i}$, or there exists $\gamma_{i}\in \mathcal{K}[g](x_{i})$, $i=1,\cdots,n$
such that $\sum_{j=1}^{n}l_{ij}\gamma_{j}=0$ for each $i$ satisfying $g$ is discontinuous at $x_{i}$.
Then let $\gamma(t)=[\gamma_{1}(t),\cdots,\gamma_{n}(t)]^{\top}$ be such that $\gamma_{i}(t)=\mathcal{K}[g](x_{i}(t))$, $i=1,\cdots,n$, and $\sum_{j=1}^{n}l_{ij}\gamma_{j}(t)=0$ for each
$i$ satisfying $g$ is discontinuous at $x_{i}(t)$, we have
\begin{eqnarray*}
\frac{dV_{L}(x(t))}{dt}&=&-\sum_{i=1}^{n}\xi_{i}[\gamma_{i}(t)-\bar{\gamma}]\sum_{j=1}^{n}l_{ij}\gamma_{j}(t)\\
&=&-\sum_{i=1}^{n}\xi_{i}[\gamma_{i}(t)-\bar{\gamma}]\sum_{j=1}^{n}l_{ij}[\gamma_{j}(t)-\bar{\gamma}]\\
&=&-\sum_{i=1}^{n}\sum_{j=1}^{n}\xi_{i}l_{ij}[\gamma_{i}(t)-\bar{\gamma}][\gamma_{j}(t)-\bar{\gamma}]\\
&=&\frac{1}{2}\Gamma(t)^{\top}(-\Xi L-L^{\top}\Xi)\Gamma(t)\\
&\le&\frac{\lambda_{2}}{2}\bar{\Gamma}(t)^{\top}\bar{\Gamma}(t),~~\text{a.e.~~}t,
\end{eqnarray*}
where $\Xi=\diag[\xi_{1},\cdots,\xi_{n}]$, $\lambda_{2}$ is the second largest eigenvalue of $-\Xi L-L^{\top}\Xi$,
$\Gamma(t)=[\gamma_{1}(t)-\bar{\gamma},\cdots,\gamma_{n}(t)-\bar{\gamma}]^{\top}$, and
$\bar{\Gamma}(t)=[\gamma_{1}(t)-\tilde{\gamma}(t),\cdots,\gamma_{n}(t)-\tilde{\gamma}(t)]^{\top}$ with $\tilde{\gamma}=\sum_{i=1}^{n}\gamma_{i}(t)/n$. The last inequality is due to the fact that largest eigenvalue of
$-\Xi L-L^{\top}\Xi$ is $0$ with the coresponding eigenspace being $k[1,1,\cdots,1]^{\top}$, $k\in\mathbb{R}$.
Since the graph of $L$ is strongly connected, $L$ is irreducible, so $\lambda_{2}<0$.
Let $\overline{i}$ be the index such that $x_{\overline{i}}(t)=\max_{i}\{x_{i}(t)\}$, and $\underline{i}$ be the index such that $x_{\underline{i}}(t)=\min_{i}\{x_{i}(t)\}$. In the case that $x_{\overline{i}}(t)>x_{\underline{i}}(t)$, we have $x_{\underline{i}}(t)<\bar{x}<x_{\overline{i}}(t)$. Otherwise, either $x_{\overline{i}}(t)\le \bar{x}$ or $x_{\underline{i}}(t)\ge \bar{x}$. If $x_{\overline{i}}(t)\le \bar{x}$, then $\Wra(x(t),L)<\bar{x}$. If  $x_{\underline{i}}(t)\ge \bar{x}$, then $\Wra(x(t),L)>\bar{x}$. These all contradict the fact that $\Wra(x(t),L)$ is constant. Thus $\gamma_{\overline{i}}(t)-\gamma_{\underline{i}}(t)>g(\bar{x}^{+})-g(\bar{x}^{-})>0$.
On the other hand,
\begin{eqnarray*}
\bar{\Gamma}(t)^{\top}\Gamma(t)&=&\sum_{i=1}^{n}[\gamma_{i}(t)-\tilde{\gamma}(t)]^{2}\\
&\ge& [\gamma_{\overline{i}}(t)-\tilde{\gamma}(t)]^{2}+[\gamma_{\underline{i}}(t)-\tilde{\gamma}(t)]^{2}\\
&\ge&\frac{1}{2}[\gamma_{\overline{i}}(t)-\gamma_{\underline{i}}(t)]^{2}\\
&>&[g(\bar{x}^{+})-g(\bar{x}^{-})]^{2}/2.
\end{eqnarray*}

Thus, we have
\begin{eqnarray*}
\frac{dV_{L}(x(t))}{dt}<-\frac{\lambda_{2}}{4}[g(\bar{x}^{+})-g(\bar{x}^{-})]^{2},~\text{a.e.}~t
\end{eqnarray*}
for $V_{L}>0$. This implies that $V_{L}$ will converge to zero in finite time upper bounded by $\displaystyle\frac{4V_{L}(x(0))}{\lambda_{2}[g(\bar{x}^{+})-g(\bar{x}^{-})]^{2}}$. The proof is completed.\qquad$\square$

\section{Proof of Theorem \ref{thmMain}}\label{appenProofThm2}
Let $V^{*}$, $V_{*}$ and $V$ be defined as in
the previous section. Given any initial value
$x(0)\in\mathbb{R}^{n}$ and any switching sequence of time points,
denoted by $0=t_{0}<t_{1}<t_{2}<\cdots$, we can construct the
solution in the following way. First, with initial value $x(0)$,
there exists a Filippov solution $x(t)$ on some interval
$[0,\delta)\subset [0,t_{1}]$. By similar arguments used in the
proof of Lemma \ref{lemShrinkingOfFilippovSolution}, we can prove
that $x(t)$ is shrinking and can be extended to the whole interval
$[0,t_{1}]$.
Repeating such arguments, we can show that %on each interval
a solution of (\ref{sysContinuousSwitchingTopology}) can be defined
as follows:
\begin{eqnarray*}
x(t)=x^{k}(t), \quad t\in [t_{k},t_{k+1}],
\end{eqnarray*}
where $x^{k}(t)$ is a Filippov solution successively defined from
$x^{k-1}(t)$ on $[t_{k},t_{k+1}]$ such that
$x^{k}(t_{k})=x^{k-1}(t_{k})$. It is obvious that $x(t)$ is
shrinking and absolutely continuous. Let $i^{*}$, $i_{*}$ be the
indices satisfying $V^{*}(x)=x_{i^{*}}$, $V_{*}(x)=x_{i_{*}}$,
respectively. Similar to the arguments in previous section, on each
interval $[t_{k},t_{k+1}]$, we have
%\begin{eqnarray*}
%\frac{dV}{dt}&=&\sum_{j=1}^{n}l^{k}_{i^{*}j}\gamma_{j}(t)-\sum_{j=1}^{n}l^{k}_{i_{*}j}\gamma_{j}(t),
%\end{eqnarray*}
%where $\gamma_{j}(t)\in\mathcal{K}[g](x_{j}(t))$ for each $j$. So we have:
\begin{eqnarray}
\frac{dV}{dt}&=&-\sum_{j=1}^{n}l^{k}_{i^{*}j}\gamma_{j}(t)+\sum_{j=1}^{n}l^{k}_{i_{*}j}\gamma_{j}(t)\nonumber\\
&=&-\sum_{j=1,j\ne i^{*}}^{n}l^{k}_{i^{*}j}[\gamma_{j}(t)-\gamma_{i^{*}}(t)]
+\sum_{j=1,j\ne i_{*}}^{n}l^{k}_{i_{*}j}[\gamma_{j}(t)-\gamma_{i_{*}}(t)]\nonumber\\
&=&(l^{k}_{i^{*}i_{*}}+l^{k}_{i_{*}i^{*}})[\gamma_{i^{*}}(t)-\gamma_{i_{*}}(t)]
+\sum_{j=1,j\ne
i^{*},i_{*}}\{l^{k}_{i^{*}j}[\gamma_{i^{*}}(t)-\gamma_{j}(t)]+l^{k}_{i_{*}j}[\gamma_{j}(t)-\gamma_{i_{*}}(t)]\}
\nonumber\\
&\le&(l^{k}_{i^{*}i_{*}}+l_{i_{*}i^{*}})[\gamma_{i^{*}}(t)-\gamma_{i_{*}}(t)]
+\sum_{j=1,j\ne i^{*},i_{*}}\max\{l^{k}_{i^{*}j},l^{k}_{i_{*}j}\}[\gamma_{i^{*}}(t)-\gamma_{i_{*}}(t)]\nonumber \\
&\le&-\epsilon(-l^{k}_{i^{*}i_{*}}-l^{k}_{i_{*}i^{*}}+\sum_{j=1,j\ne
i^{*},i_{*}}\min\{-l^{k}_{i^{*}j},-l^{k}_{i_{*}j}\})V\nonumber\\
&\le&-\epsilon\eta(-L^{k})V,
\end{eqnarray}
where $\gamma_{j}(t)\in\mathcal{K}[g](x_{j}(t))$ for each $j$.
Therefore, we have
\begin{eqnarray}
V(x(t_{k+1}))\le e^{-\epsilon\eta(-L^{k})\Delta t_{k}}V(x(t_{k}))\le
e^{-\epsilon\sum_{i=0}^{k}\eta(-L^{i})\Delta t_{i}}V(x(0)). \label{sto}
\end{eqnarray}
Thus if
$\sum_{k=1}^{+\infty}\eta(-L^{k})\Delta t_{k}=\infty$, then
\begin{align}
\lim_{k\to+\infty}V(x(t_{k}))=0.
\end{align}

On the other hand, let $\mathbb{S}_{\mathbb{N}}$ denote the space of strictly increasing infinite sequence of
the natural numbers, we have
\begin{eqnarray}
\Prb\big\{\sum_{k=1}^{+\infty}\eta(-L^{k})\Delta t_{k}=\infty\big\}
&\ge& \Prb\big\{\eta(-L^{n_{k}})\ge \delta,\sum_{k=1}^{\infty}\Delta t_{n_{k}}=\infty, \{n_{k}\}\in \mathbb{S}_{\mathbb{N}}\big\}\nonumber\\
&=&\sum_{\{n_{k}\}\in \mathbb{S}_{\mathbb{N}}}\Prb\{\eta(-L^{n_{k}})\ge \delta\big|\sum_{k=1}^{\infty}\Delta t_{n_{k}}=\infty\}\Prb\{\sum_{k=1}^{\infty}\Delta t_{n_{k}}=\infty\}\label{eqn1}\\
&=&\sum_{\{n_{k}\}\in \mathbb{S}_{\mathbb{N}}}\Prb\{\eta(-L^{n_{k}})\ge \delta\}\Prb\{\sum_{k=1}^{\infty}\Delta t_{n_{k}}=\infty\}\label{eqn2}\\
&=&\sum_{\{n_{k}\}\in \mathbb{S}_{\mathbb{N}}}\Prb\{\eta(-L^{n_{k}})\ge \delta\}\label{eqn3}\\
&=&\Prb\{\eta(-L^{k})\ge \delta \text{ infinitely often}\}\nonumber\\
&=&1\nonumber.
\end{eqnarray}
Due to the independence of $\{\Delta t_{i}\}$ and
$\{L^{k}\}$ from Assumption \ref{asumpIndependentConditions}, we can have the equality from \eqref{eqn1} to
\eqref{eqn2}. Since $\{\Delta t_{k}\}$ is independent and identically distributed, the subsequence $\{\Delta t_{n_{k}}\}$ is also independent and identically distributed for each $\{n_{k}\}\in \mathbb{S}_{\mathbb{N}}$.
From the strong law of large numbers, we have
$$
\lim_{N\to \infty}\frac{1}{N}\sum_{k=1}^{N}\Delta t_{n_{k}}=\Exp\Delta t_{1}>0
$$
almost surely, which implies
$$
\Prb\{\sum_{k=1}^{\infty}\Delta t_{n_{k}}=\infty\}=1.
$$

Thus we get the equality from \eqref{eqn2} to
\eqref{eqn3}.

This implies
\begin{eqnarray*}
P\{\sum_{k=1}^{+\infty}\eta(-L^{k})\Delta t_{k}=\infty\}=1,
\end{eqnarray*}
and
\begin{eqnarray*}
P\{\lim_{k\to+\infty}V(x(t_{k}))=0\}=1.
\end{eqnarray*}
because $V(x(t))$ is nonincreasing with respect to $t$, we conclude
\begin{eqnarray*}
P\{\lim_{t\to+\infty}V(x(t))=0\}=1,
\end{eqnarray*}
Theorem \ref{thmMain} is proved completely.\qquad$\square$

\end{document}